\documentclass [11pt,twoside,a4paper]{article}
\usepackage{amsfonts}
\usepackage{amsthm}
\usepackage{amsmath}
\usepackage{amstext}
\usepackage{amssymb}
\usepackage{mathrsfs}
\usepackage{amscd}
\usepackage{xypic}
\usepackage{epsf}              
\usepackage{graphicx}          
\usepackage{fancybox}          
\usepackage{color}             
\usepackage{fancyhdr}
\usepackage[hang,footnotesize]{caption2}  

\def\mathcal{\mathscr}
\newfont{\aaa}{cmb10 at 19pt}
\newfont{\bbb}{cmb10 at 11pt}
\newtheorem{lemma}{Lemma}[section]
\newtheorem{theorem}{Theorem}[section]
\newtheorem{corollary}{Corollary}[section]

\newcommand{\beq}{\begin{equation}}
\newcommand{\eeq}{\end{equation}}
\newcommand{\bey}{\begin{eqnarray}}
\newcommand{\eey}{\end{eqnarray}}
\newcommand{\beyy}{\begin{eqnarray*}}
\newcommand{\eeyy}{\end{eqnarray*}}

\makeatletter
\def\@evenhead{
\vbox{\hbox to \textwidth {}{\hspace{0mm}{\footnotesize
\thepage}}{\hspace{8cm} {\footnotesize {Chuancun YIN, Kam C. YUEN}}} \protect\vspace{1truemm}\relax \hrule depth0pt
height0.15truemm width\textwidth}}
\def\@evenfoot{}
\def\@oddhead{\vbox{\hbox to \textwidth
{{\hspace{0cm}{\footnotesize  Exact joint laws associated with spectrally
 negative L\'evy processes }\hfill{\footnotesize
\thepage}}\hspace{0mm}}{} \protect\vspace{1truemm}\relax\hrule
depth0pt height0.15truemm width\textwidth}}
\def\@oddfoot{}
\makeatother

\begin{document}

\thispagestyle{empty} \thispagestyle{fancy} {
\fancyhead[lO,RE]{\footnotesize  Front. Math. China \\
DOI 10.1007/s11464-012-0186-5\\[3mm]
\includegraphics[0,-50][0,0]{11.bmp}}
\fancyhead[RO,LE]{\scriptsize \bf 
} \fancyfoot[CE,CO]{}}
\renewcommand{\headrulewidth}{0pt}


\setcounter{page}{1}
\qquad\\[8mm]

\noindent{\aaa{Exact joint laws associated with spectrally
 negative L\'evy processes and applications to insurance risk theory}}\\[1mm]

\noindent{\bbb Chuancun YIN}\\[-1mm]

\noindent\footnotesize{ School of Mathematical Sciences, Qufu Normal University\\
         Shandong 273165, China}\\[6mm]
\noindent{\bbb   Kam C. YUEN}\\[-1mm]

\noindent\footnotesize{ Department of Statistics and Actuarial Science\\
The University of Hong Kong,
Pokfulam, Hong Kong, China}\\[6mm]

\vskip-2mm \noindent{\footnotesize$\copyright$ Higher Education
Press and Springer-Verlag Berlin Heidelberg 2013} \vskip 4mm

\normalsize\noindent{\bbb Abstract}\quad We consider the  spectrally negative L\'evy processes and
determine the joint laws for the quantities such as the first and
last passage times over a fixed level, the overshoots and
undershoots at first passage, the minimum, the maximum and the
duration of negative values. We apply our  results to insurance risk
theory to find an explicit expression for the generalized expected
discounted penalty function in terms of scale functions. Further, a
new expression for the generalized Dickson's formula is provided.\vspace{0.3cm}

\footnotetext{Received May 10, 2012; accepted October 10, 2013\\
\hspace*{5.8mm}Corresponding author: Chuancun YIN, E-mail:
ccyin@mail.qfnu.edu.cn}

\noindent{\bbb Keywords}\quad Fluctuation
identity,  spectrally negative L\'evy processes, generalized
Dickson's formula,  scale functions,   occupation times, Suprema and
infima\\
{\bbb MSC}\quad 60G51; 60G50; 60J75; 91B30\\[0.4cm]

\section {Introduction}
Let $X=\{X(t), t\ge 0\}$ be a real valued spectrally negative L\'evy
process, i.e. a stochastic process with c\`{a}dl\`{a}g paths without
positive jumps that has stationary independent increments defined on
some filtered space $(\Omega, \cal{F}$, $\{{{\cal F}_t}, t\ge 0\}$, $P)$
where the filtration $\{{\cal{F}}_t, t\ge 0\}$  satisfies the usual
conditions of right continuity and completion. The reader is
referred to Bertoin [5],  Kyprianou [25] and Doney [12] for a
general discussion to L\'evy processes. Write $P_x$ for the
probability law of $X$ when $X(0)=x$ and $E_x$ for the expectation
with respect to $P_x$. We simply write $P_0=P$ and $E_0=E$. As
usual, we will exclude the case that $X$ has monotone paths. By the
absence of positive jumps, the moment generating function of $X(t)$
exists for all $\alpha\ge 0$ and is given by $E\text{e}^{\alpha
X(t)}=\text{e}^{t\varphi(\alpha)}, \; t\ge 0,$ where $\varphi:
[0,\infty)\rightarrow\Bbb{R}$
 is known as the Laplace exponent. It is given by the L\'evy-Khinchin formula
\begin{equation}
\varphi(\alpha)=\mu \alpha+\frac12\sigma^2\alpha^2+\int_{-\infty}^0(\text{e}^{\alpha
x}-1-\alpha x 1_{(x>-1)})\Pi(dx),
\end{equation}
  where $\mu\in\Bbb{R}$,
$\sigma\ge 0$ and $\Pi$ is a measure on $(-\infty,0)$ known as the
L\'evy measure satisfying $\int_{-\infty}^0 (x^2\wedge
1)\Pi(dx)<\infty$. It is known that $\varphi(x)$ is strictly convex
for $x\ge 0$ and that $\varphi(0)=0$ and
$\lim_{x\to\infty}\varphi(x)=\infty$. Moreover,
$\varphi'(0+)=E(X(1))=\mu+\int_{-\infty}^{-1}x\Pi(dx)\in
[-\infty,\infty)$. Recall that $X$ drifts to $-\infty$ if and only
if $\varphi'(0+)<0$. Such a spectrally negative L\'evy process has
sample paths of bounded variation if and only if $\sigma=0$ and
$\int_{-1}^0|x|\Pi(dx)<\infty.$

Define the right inverse $\Phi(q)=\sup\{\alpha\ge 0:
\varphi(\alpha)=q\}$ for each $q\ge 0$. If $\varphi'(0+)\ge 0$ then
$\alpha=0$ is the unique solution to $\varphi(\alpha)=0$ and
otherwise there are two solutions with $\alpha=\Phi(0)>0$ being the
larger of the two, the other is $\alpha=0$. The Laplace exponent
$\varphi$ is continuous and increasing over $[\Phi(0),\infty)$, so
that $\lim_{\alpha\to 0+}\Phi(\alpha)=\Phi(0).$

It is  known that (see Ref. [3]) for any $c$ such that
$\varphi(c)=\log E[\exp cX(1)]$ is finite, the family
$\Lambda_t(c)=\exp\{c(X(t)-x)-\varphi(c)t\}$ is a martingale under
$P_x$. Let $P^{(c)}_x$ denote the probability measure on
$\cal{F}_t=\sigma$ $\{X(s):0\le s<\infty\}$ defined by
\begin{equation}
\frac{\text{\rm d} P^{(c)}_x}{\text{\rm d}P_x}\bigg{\vert}_{{\cal F}_t}=\frac{\Lambda_t(c)}{\Lambda_0(c)}
\end{equation}
for all $0\le t<\infty$.   Under the measure $P^{(c)}_x$, $X$
remains within the class of spectrally negative process and the
Laplace exponent of $X$ is given by
$\varphi_c(\alpha)=\varphi(\alpha+c)-\varphi(c),\ \alpha\ge -c$.

When we set $c=\Phi(q)$ for $q\ge 0$ we find that
$\varphi'_{\Phi(q)}(0)=\varphi'(\Phi(q))\ge 0$. Thus,
$E^{(\Phi(q))}(X(1))=\varphi'(\Phi(q))>0$ for $q>0$.  In particular,
under the measure $P^{(\Phi(q))}$, $X$ always drifts to $\infty$ for
$q>0$ (see  Kyprianou [25, pp 213-214] for more details).

Denote by $I$ and $S$ the past infimum and supremum of $X$
respectively, that is, $I_t=\inf_{0\le s\le t}X(s) \;\; \text{\rm
and}\;\; S_t=\sup_{0\le s\le t}X(s).$ Define the first passage times
above and below $a$ for $X$ by $T_a^-=\inf\{t\ge 0: X(t)< a\}\;\;
\text{\rm and}\;\;
 T_a^+=\inf\{t\ge 0: X(t)>a\}.$
We simply write $T$ for $T_0^-$.  Finally, let $T_0$ denote the time
of recovery: $T_0=\inf\{t:t>T, X(t)=0\}.$

The fluctuation theory for spectrally negative L\'evy processes has
been the object of several studies over the last 40 years, among
many others, see Emery [15], Bingham [7] and Bertoin [6]; We
refer to Kyprianou and Palmowaki [26] and Pistorius [32] for
exhaustive reviews. In a  great variety of problems of applied
fields such as actuarial mathematics, mathematical finance and
queueing theory and so on one faces the consideration of the
fluctuation theory for a class of spectrally negative L\'evy
processes; For example,  The exit problems have recently been used
by Avram et al. [3] in the context of finance in connection with
American and Canadian options. Alili and Kyprianou [1] provided,
with the help of a fluctuation identity, a generic link between a
number of known identities for the first passage time and overshoot
above/below a fixed level of a L\'evy process to the American
perpetual put optimal stopping problem. In the theory of actuarial
mathematics, the first passage problem is very important to the ruin
problems; see, for example, among others, Yang and Zhang [35],
Kl\"{u}ppelberg et al. [23], Huzak et al. [20], Zhou
[36], Garrido and Morales [17], Biffis and Kyprianou [8] and
Billis and Morales [9].

Several recent papers have concerned  the joint laws of overshoots
and undershoots of  L\'evy processes at the first and the last
passage times of a constant barrier. For example, Doney and
Kyprianou [13] studied the problem for general L\'evy process $X$.
They computed the quintuple law of the time of the first passage
relative to the time of the last maximum at the first passage, the
time of the last maximum at the first passage, the overshoot at the
first passage, the undershoot at the first passage and the
undershoot of the last maximum at the first passage, i.e. the law of
$(T_x^+ -\overline{G}_{{T_x^+}-}, \overline{G}_{{T_x^+}-},
X_{T_x^+}-x, x-X_{{T_x^+}-}, x-S_{{T_x^+}-}),$ where
$\overline{G}_t=\sup\{s<t, X_s=S_s\}$. Kyprianou, Pardo and Rivero
(2010) extended above quintuple law to a family of related joint
laws for L\'evy processes, L\'evy processes conditioned to stay
positive and positive self-similar Markov processes at both first
and last passage over a fixed level. Eder and  Kl\"uppelberg [16]
derived new results in fluctuation theory for sums of possibly
dependent L\'evy processes.   Chaumont, Kyprianou and Pardo [10]
consider some special classes of L\'evy processes with no Gaussian
component whose L\'evy measure is of the type $\Pi(dx)=\pi(x)dx$,
where $\pi(x)=e^{\gamma x}\nu (e^x-1)$.

Motivated by these interesting papers, we continue to study the
fluctuation theory of  L\'evy processes. In this paper, we restrict
ourselves to the spectrally negative L\'evy processes. The advantage
is that  all results can be explicitly expressed in terms of scale
functions. The rest of the paper is organized as follows. The next
section  reviews some preliminary results of spectrally negative
L\'evy processes that will be needed later on. In Section 3 we
determine the joint laws of all or some of the quantities such as
the first and last passage time over a fixed level, the overshoots
and undershoots at the first passage, the minimum, the maximum and
the duration of negative values. Applications in insurance risk
theory are discussed in Section 4.

\section  {Preliminaries}
\setcounter{equation}{0}

 Let us now reviewing some preliminary results of spectrally negative
L\'evy processes. We will assume that the measure $\Pi$ has a
density $\pi$, with respect to the Lebesgue measure. So that the
scale functions are differentiable. For details, see  Ref. [29].\\

{\bf 2.1 The scale functions and the survival probabilities}\\

Scale functions are  key object in the fluctuation theory of
spectrally negative L\'evy processes and its applications, the
survival probabilities are  key object in risk theory.

{\bf Definition 2.1} For $q\ge 0$, the $q$-scale function $W^{(q)} :
(-\infty,\infty)\rightarrow [0,\infty)$ is the unique function whose
restriction to $(0,\infty)$ is continuous and has Laplace transform
$$\int_0^{\infty}e^{-\alpha x}W^{(q)}(x)dx =\frac{1}{\varphi(\alpha)-q},\;\alpha>\Phi(q)$$
and is defined to be identically zero for $x< 0$. For short we shall
write $W^{(0)}=W$. Further, we shall use the notation $W^{(q)}_c(x)$
to mean the $q$-scale function as defined above for $(X, P^{(c)})$.

 For $q\ge 0$, we define $Z^{(q)}(x)=1$ for
$x\le 0$ and
$$Z^{(q)}(x)=1+q\int^x_0 W^{(q)}(y)dy \;\text{for}\; x>0.$$

For $x\ge 0$, define the survival probability
$$\overline{Q^{(c)}}(x)=1-Q^{(c)}(x)=P_x^{(c)}(I_{\infty}\ge 0),$$ where
$\overline{Q^{(c)}}(0)=\lim_{x\downarrow 0}\overline{Q^{(c)}}(x)$.
When $c=0$, we write $\overline{Q}(x)$ instead of
$\overline{Q^{(0)}}(x)$.

It is well known that (see  Ref. [26]), under the
condition  $\varphi'(0+)>0$, $W(x)=\frac{1}{\varphi'(0+)}
P_x(I_{\infty}\ge 0)$.  For  $q>0$ ($\varphi'(0+)>0$ is not
required),  using the fact that
$W^{(q)}(x)=\text{e}^{\Phi(q)x}W_{\Phi(q)}(x)$, we also have
\begin{equation}
W^{(q)}(x)=\frac{1}{\varphi'(\Phi(q))}\text{e}^{\Phi(q)x}\overline{Q^{(\Phi(q))}}(x).
\end{equation}\\

{\bf 2.2 The triple law of $T$, $X(T-)$ and $X(T)$}\\

We denote the $q$-potential measure of $X$ killed on exiting
$[0,\infty)$ with starting point $x$ by $U^{(q)}(x,dy)$. That is
$U^{(q)}(x,dy)=dy\int_0^{\infty}\text{e}^{-q t}P_t(x,y)dt,$ for
$q\ge 0$ with $U^{(0)}=U$, where $P_t(x,y)dy=P_x(T>t, X(t)\in dy).$
 If a density of $U^{(q)}(x,dy)$ exists
with respect to the Lebesgue measure for each $x\ge 0$ then we call
it the potential density and denoted it by $u^{(q)}(x,y)$ (with
$u^{(0)}=u)$. It is well known that (see Refs. [5, 12, 25])
 \begin{equation}
 u^{(q)}(x,y)=W^{(q)}(x)\text{e}^{-\Phi(q) y}-W^{(q)}(x-y).
\end{equation}
For $X(0)=x\ge 0$, let
 \begin{equation}
 f(y,z,t|x)dydzdt=P_x( T\in dt, X(T-)\in dy, |X(T)|\in dz).
 \end{equation}
For $q\ge 0$, define
$$f_q(y,z|x)=\int^{\infty}_0\text{e}^{-q t}f(y,z,t|x)\text{\rm d}t,\;\;
f_q(y|x)=\int_0^{\infty}f_q (y,z|x)\text{\rm d}z.$$

It follows from Doney [12, Remark 5(i), p105] that (by letting
$a\to\infty$)  for $x, y, z>0$,
\begin{equation}
E_x(\text{e}^{-q T}, X(T-)\in dy, |X(T)|\in dz)=u^{(q)}
(x,y)\pi(-z-y)dzdy,
\end{equation}
and
$$\int_0^{\infty}\text{e}^{-q
t}f(y,z,t|x)dt=\pi(-z-y)\int_0^{\infty}\text{e}^{-q t}P_t(x,y)dt.$$
So that
\begin{equation}
f_q(y,z|x)=u^{(q)}(x,y)\pi(-z-y),
\end{equation}
and
$$f(y,z,t|x)=P_t(x,y)\pi(-z-y).$$

From (2.1), (2.2)  and (2.5) we obtain for $y>0, x>0$,
 \begin{eqnarray}
 f_q(y|x)&=&(\Pi(-\infty,-y))(W^{(q)}(x)e^{-\Phi(q)
y}-W^{(q)}(x-y))\nonumber\\
&=&(\Pi(-\infty,-y))\frac{e^{\Phi(q)(x-y)}}
{\varphi'(\Phi(q))}\left(\overline{Q^{(\Phi(q)}}(x)-\overline{Q^{(\Phi(q)}}(x-y)
\right).
\end{eqnarray}
  If the paths of $X$ are of  bounded
variation,  then by (2.5) and (2.6)
\begin{equation}
f_q(y,z|0)=b^{-1}\text{e}^{-\Phi(q) y}\pi(-z-y),\;\; f_q(y|0)=(\Pi(-\infty, -y))
b^{-1}\text{e}^{-\Phi(q) y}
\end{equation}
since $W^{(q)}(0+)=b^{-1}$,
where $b=\mu-\int_{-1}^0 x\Pi(dx)$. It follows from (2.5)-(2.7) that
\begin{equation}\frac{f_{q}(y,z|x)}{f_{q}(y,z|0)}=\frac{f_{q}(y|x)}{f_{q}(y|0)}=
  \left\{
  \begin{array}{ll}
      bW^{(q)}(x),&0\le x\le y,\\
     b(W^{(q)}(x)-e^{\Phi(q) y}W^{(q)}(x-y)), & x\ge y>0.
  \end{array}
  \right.
\end{equation}
 Finally, for $q\ge 0$, define $\Phi_c(q)$ to be the largest
real root of the equation $\varphi_c(\theta)=q$. Then we have the
following important result which is due to Emery [15]; See also Kyprianou [25].
\begin{lemma}
For any $\alpha\ge 0$ and $\beta\ge 0$, the
joint Laplace transform of $T_y^-$ and $X(T_y^-)$, with the initial
condition $X(0)=x>y$,  is given by
\begin{equation}
E_x(\text{e}^{-\alpha T^-_y+\beta X(T_y^-)}, T<\infty)
=e^{\beta
x}\left(Z_{\beta}^{(p)}(x-y)-\frac{p}{\Phi_{\beta}(p)}W_{\beta}^{(p)}(x-y)\right),
\end{equation}
where $W_{\beta}^{(p)}$ and $Z_{\beta}^{(p)}$ are scale functions
with respect to the measure $P^{(\beta)}, p=\alpha-\varphi(\beta)$, $\Phi_{\beta}(p)=\Phi(\alpha)-\beta$
and $\frac{p}{\Phi_{\beta}(p)}$ is understand in the limiting sense
if $p=0$.
\end{lemma}

\section  {Joint laws for the spectrally negative L\'evy processes}
\setcounter{equation}{0}

The main purpose of this section is to investigate  some joint laws
for the spectrally negative L\'evy process involving some or all of
the first passage time, the last passage time, the overshoots and
undershoots at first passage, the minimum, the maximum and the
duration of negative values.

\begin{theorem}   For $q\ge 0$ and  for positive numbers
$x,y, z, a$ and $b$ such that $b<a\wedge x\wedge y$, $a>x, a>y$,
  \begin{eqnarray}
   &E_x&(\text{e}^{-q T}, X(T-)\in dy, |X(T)|\in dz, I_{T-}>b,
S_{T-}\le a,
T<\infty)\nonumber\\
&=&\pi(-z-y)\left(\frac{W^{(q)}(x-b)}{W^{(q)}(a-b)}
W^{(q)}(a-y)-W^{(q)}(x-y)\right)dydz.
\end{eqnarray}
\end{theorem}
{\bf Proof} Using the spatial homogeneity and the strong Markov
property of $X(t)$,  we obtain
 \begin{eqnarray}
 &&P_x(X(T-)\in dy, |X(T)|\in dz,  I_{T-}>b,  S_{T-}>a, T<\infty)\nonumber\\
&=&P_{x-b}(X(T-)\in dy-b, |X(T)|\in dz+b, I_{T-}>0, X(T)\le -b,S_{T-}>a-b,  T<\infty)\nonumber\\
&=&P_{x-b}(X(T-)\in dy-b, |X(T)|\in dz+b, S_{T-}>a-b,  T<\infty)\nonumber\\
&&-P_{x-b}(X(T-)\in dy-b, |X(T)|\in dz+b, X(T)> -b,S_{T-}>a-b,  T<\infty)\nonumber\\
&=&P_{x-b}(S_{T-}>a-b)P_{a-b}(X(T-)\in dy-b, |X(T)|\in dz+b,   T<\infty)\nonumber\\
&=&\frac{W(x-b)}{W(a-b)}f_0(y-b,z+b|a-b)dydz, \nonumber
\end{eqnarray}
 \begin{eqnarray}
  &&P_x(X(T-)\in dy, |X(T)|\in dz, I_{T-}>b, T<\infty)\nonumber\\
&=&P_{x-b}(X(T-)\in dy-b, |X(T)|\in dz+b, I_{T-}>0, X(T)\le -b,  T<\infty)\nonumber\\
&=&P_{x-b}(X(T-)\in dy-b, |X(T)|\in dz+b,   T<\infty)\nonumber\\
&&-P_{x-b}(X(T-)\in dy-b, |X(T)|\in dz+b,  X(T)>-b,  T<\infty)\nonumber\\
&=&f_0(y-b,z+b|x-b)dydz,\nonumber
\end{eqnarray}
  where
\begin{eqnarray*}
  f_0(y,z|a)dydz&=&P_a(X(T-)\in dy, |X(T)|\in dz, T<\infty)\\
&=&\pi(-z-y)(\text{e}^{-\Phi(0) y}W(a)-W(a-y))dydz.
\end{eqnarray*}
It follows that
 \begin{eqnarray}
  &P_x&(X(T-)\in dy, |X(T)|\in dz, I_{T-}>b, S_{T-}\le a,
T<\infty)\nonumber\\
&=&\pi(-z-y)\left(\frac{W(x-b)}{W(a-b)}W(a-y)-W(x-y)\right)dydz.
\end{eqnarray}
To prove (3.1), using (3.2) and applying the exponential change of
measure   we get
 \begin{eqnarray*}
  &E_x&(\text{e}^{-q T}, X(T-)\in dy, |X(T)|\in dz, I_{T-}>b,
S_{T-}\le a,
T<\infty)\\
&&=\Pi_{\Phi(q)}(-dz-y)\text{e}^{-\Phi(q)(-z-x)}\\
&&\times\left(\frac{W_{\Phi(q)}(x-b)}{W_{\Phi(q)}(a-b)}W_{\Phi(q)}(a-y)-W_{\Phi(q)}(x-y)\right)dy,
\end{eqnarray*}
where $W_{\Phi(q)}$ is the 0-scale function of $X$ under
$P^{(\Phi(q))}$. Recall that
$$W^{(q)}(x)=\text{e}^{\Phi(q)x}W_{\Phi(q)}(x),\;\; \Pi_{\Phi(q)}(-dx)=\text{e}^{-\Phi(q)x}\Pi(-dx),$$
where $\Pi_{\Phi(q)}$ is the L\'evy measure of $U$ under
$P^{(\Phi(q))}$, we complete the proof.

Letting $a\to\infty$ or $b\to 0$ in (3.1) and using the  fact in
Zhou [37]
$$\lim_{a\to\infty}\frac{W^{(q)}(a-y)}{W^{(q)}(a-b)}=e^{-(y-b)\Phi(q)}$$
 yield the following result.

\begin{corollary}  For $q\ge 0$ and  for positive numbers
$x,y,z, a$ and $b$ such that $b< a\wedge x\wedge y$, $a>x, a>y$,
  \begin{eqnarray}
   (1).\;\; &&E_x(\text{e}^{-q T}, X(T-)\in dy, |X(T)|\in dz, I_{T-}>b,
T<\infty)\nonumber\\
&&=\pi(-z-y)\left(W^{(q)}(x-b)e^{-(y-b)\Phi(q)}-W^{(q)}(x-y)\right)dydz.
 \end{eqnarray}
 \begin{eqnarray} (2).\;\; &&E_x(\text{e}^{-q T}, X(T-)\in dy, |X(T)|\in dz, S_{T-}\le a,
T<\infty)\nonumber\\
&&=\pi(-z-y)\left(\frac{W^{(q)}(x)}
{W^{(q)}(a)}W^{(q)}(a-y)-W^{(q)}(x-y)\right)dydz.
 \end{eqnarray}
\end{corollary}
 {\bf Remark 3.1 } Taking derivative with respect to $b$ in
 (3.3) yields the following result:
 \begin{eqnarray}
 &&E_x(e^{-q T}, X(T-)\in dy, |X(T)|\in dz, I_{T-}\in db,
T<\infty)\nonumber\\
&&=\pi(-z-y)e^{-\Phi(q)(y-b)}({W^{(q)}}'(x-b)-\Phi(q)
W^{(q)}(x-b))dydzdb. \nonumber
\end{eqnarray}
 In the case that $X$
drifts to $\infty$,  Biffis and Kyprianou [8] also found the
result based on a quintuple law in Doney and Kyprianou [13].

\begin{theorem}   For $q, \beta\ge 0$ and  for positive
numbers $x,y, z, a$ and $b$ such that $z<b\le a\wedge x\wedge y$,
$a>x, a>y$,
\begin{eqnarray} &&E_x(\text{e}^{-q T-\beta (T_0-T)}, X(T-)\in dy, |X(T)|\in dz, I_{T_0}>-b,
S_{T_0}\le a, T_0<\infty)\nonumber\\
&&=K_x(y,z,a) \frac{W^{(\beta)}(-z+b)}{W^{(\beta)}(b)}dydz,
\end{eqnarray}
where
$$K_x(y,z,a)=\pi(-z-y)\left(\frac{W^{(q)}(x)}{W^{(q)}(a)}
W^{(q)}(a-y)-W^{(q)}(x-y)\right).$$
\end{theorem}
 {\bf Proof} Applying  the strong Markov property of $\{X(t)\}$ at
$T$ we get that
 \begin{eqnarray}
  &&E_x(\text{e}^{-q T-\beta (T_0-T)}, X(T-)\in dy, |X(T)|\in dz, I_{T_0}>-b,
S_{T_0}\le a, T_0<\infty)\nonumber\\
&&=E_x(\text{e}^{-q T-\beta (T_0-T)}, X(T-)\in dy, |X(T)|\in dz,\nonumber\\
&&\inf_{T\le t\le T_0}X(t)>-b, S_{T-}\le a, T_0<\infty)\nonumber\\
&&=E_x(\text{e}^{-q T}, X(T-)\in dy, |X(T)|\in dz, S_{T-}\le a,
T<\infty)\nonumber\\
&&\times E_{-z}(\text{e}^{-\beta T_0}, I_{T_0}>-b,
T_0<\infty).
 \end{eqnarray}
Applying the
exponential change of measure  one gets
 \begin{eqnarray} E_{-z}(\text{e}^{-\beta T_0},
I_{T_0}>-b, T_0<\infty)
&=&\text{e}^{-\Phi(\beta)z}E_{-z}^{\Phi(\beta)}(I_{T_0}>-b, T_0<\infty)\nonumber\\
&=&\text{e}^{-\Phi(\beta)z}
\frac{W_{\Phi(\beta)}(-z+b)-W_{\Phi(\beta)}(-z)}{W_{\Phi(\beta)}(b)}\nonumber\\
&=&\frac{W^{(\beta)}(-z+b)}{W^{(\beta)}(b)}.
 \end{eqnarray}
 Now (3.5) follows from (3.4), (3.6) and (3.7). This completes the proof of Theorem 3.2.

Letting $a\to\infty$ and $b\to \infty$ in (3.5)
 yield the following result.

\begin{corollary}  For $q,\beta\ge 0$ and  for positive
numbers $x,y,z$  we have
 \begin{eqnarray}
 &&E_x(\text{e}^{-q T-\beta(T_0-T)}, X(T-)\in dy, |X(T)|\in dz,
T_0<\infty)\nonumber\\
&&=\pi(-z-y)e^{-\Phi(\beta)z}u^{(q)}(x,y)dydz.
\end{eqnarray}
\end{corollary}

The following result is an extension of Chiu and Yin [11, Theorem 3.2].
 \begin{theorem}
 Suppose that the L\'evy process $X$
with the Laplace exponent (1.1) drifts to $\infty$. Denote by
$l=\sup\{t\ge 0: X(t)<0\}$ the last passage time below level 0.  For
$q, \beta\ge 0$ and for positive numbers $x,y, z, a$ and $b$ such
that $z<b\le a\wedge x\wedge y$, $a>x, a>y$, then
 \begin{eqnarray}
 &&E_x(\text{e}^{-q T-\beta (l-T)}, X(T-)\in dy, |X(T)|\in dz, I_{l}>-b,
S_{l}<a, T<\infty)\nonumber\\
&&=K_x(y,z,a) R(z,a,b)dydz,
\end{eqnarray}
 where
$$K_x(y,z,a)=\pi(-z-y)\left(\frac{W^{(q)}(x)}{W^{(q)}(a)}
W^{(q)}(a-y)-W^{(q)}(x-y)\right),$$
  \begin{eqnarray*}
   R(z,a,b)&=&\varphi'(0+)e^{-\Phi(\beta)(b-z)}W^{(\beta)}(b-z)-\varphi'(0+)\Phi'(\beta)\\
&&+\varphi'(0+) \frac{W^{(\beta)}(b-z)}{W^{(\beta)}(a+b)}\left(W^{(\beta)}(a)+\Phi'(\beta)-e^{-(a+b)\Phi(\beta)}
W^{(\beta)}(a+b)\right).
\end{eqnarray*}
 \end{theorem}
{\bf Proof} \ The strong Markov property of $X$ yields that the left
hand side of (3.9) is equals to
 \begin{eqnarray*}
 &&E_x(\text{e}^{-q T},
X(T-)\in dy, |X(T)|\in dz, S_{T}< a, T<\infty)\\
&&\times E_{-z}(\text{e}^{-\beta l},  I_{l}>-b, S_{l}<a,
l>0):=S_x(y,z,a)\times R(z,a,b).
\end{eqnarray*}
 Letting $\beta\to 0$ and
$b\to\infty$ in (3.5) yield $S_x(y,z,a)=K_x(y,z,a)$ since $K_x(y,z,\{a\})=0$. Applying the
strong Markov property of $X$  one finds
 \begin{eqnarray} R(z,a,b)&=&E_{-z}(\text{e}^{-\beta l}, l<T^+_a,T^-_{-b}=\infty)\nonumber\\
&=&E_{-z}(\text{e}^{-\beta l}, T^-_{-b}=\infty)-E_{-z}(\text{e}^{-\beta l}, l>T^+_a,T^-_{-b}=\infty)\nonumber\\
&=&E_{-z}(\text{e}^{-\beta l},l>0)-E_{-z}(\text{e}^{-\beta l}, T^-_{-b}<\infty)\nonumber\\
&&-E_{-z}\left(\text{e}^{-\beta T^+_a}E_a(\text{e}^{-\beta l},l>0,T^-_{-b}=\infty),T_a^+<T_{-b}^-\right)\nonumber\\
&=&E_{-z}(\text{e}^{-\beta l},l>0)-E_{-z}(\text{e}^{-\beta T^-_{-b}} E_{X(T^-_{-b})}(\text{e}^{-\beta l}, l>0), T^-_{-b}<\infty)\nonumber\\
&&-E_{-z}(\text{e}^{-\beta T^+_a}, T_a^+< T_{-b}^-)E_a(\text{e}^{-\beta l},
l>0,T^-_{-b}=\infty).
\end{eqnarray}
Note that
\begin{eqnarray*}
&E_a&(\text{e}^{-\beta l},l>0,T^-_{-b}=\infty)=E_a(\text{e}^{-\beta
l},l>0)\\
&&-E_a\left(\text{e}^{-\beta T^-_{-b}}E_{X(T^-_{-b})}(\text{e}^{-\beta
l},l>0), T^-_{-b}<\infty\right)\\
&=&E_a(\text{e}^{-\beta
l},l>0)-\varphi'(0+)\Phi'(\beta)E_{a}\left(\text{e}^{-\beta
T^-_{-b}+\Phi(\beta)X(T^-_{-b})}, T^-_{-b}<\infty\right)\\
&=&E_a(\text{e}^{-\beta
l},l>0)-\varphi'(0+)\Phi'(\beta)P_{a}^{(\Phi(\beta))}(T^-_{-b}<\infty),
\end{eqnarray*}
and
$$E_{-z}(\text{e}^{-\beta T^-_{-b}} E_{X(T^-_{-b})}(\text{e}^{-\beta l}, l>0),
T^-_{-b}<\infty)=\varphi'(0+)\Phi'(\beta)P_{-z}^{(\Phi(\beta))}(T^-_{-b}<\infty),$$
where we have used the change of measure argument and the formula
$$E_{u}(\text{e}^{-\beta l},l>0)=\varphi'(0+)\Phi'(\beta)\text{e}^{\Phi(\beta)u}, u\le 0.$$
See, Kyprianou [25, Ex. 8.10] or Chiu and Yin [11, (1.5) and
Theorem 3.1].

It follows that
 \begin{eqnarray}
  R(z,a,b)&=&E_{-z}(\text{e}^{-\beta l},l>0)- \varphi'(0+)\Phi'(\beta)P_{-z}^{(\Phi(\beta))}(T^-_{-b}<\infty)\nonumber\\
&&- E_{-z}(\text{e}^{-\beta T^+_a}, T_a^+< T_{-b}^-)\nonumber\\
&&\times\left(E_a(\text{e}^{-\beta
l},l>0)-\varphi'(0+)\Phi'(\beta)P_{a}^{(\Phi(\beta))}(T^-_{-b}<\infty)\right).
\end{eqnarray}

Note that (cf. Kyprianou [25, Ex. 8.10]),
 \begin{equation}
 E_{-z}(\text{e}^{-\beta l},l>0)=\varphi'(0+)\Phi'(\beta)\text{e}^{-\Phi(\beta)z},
 \end{equation}
 \begin{equation}
E_a(\text{e}^{-\beta
l},l>0)=\varphi'(0+)\Phi'(\beta)\text{e}^{a\Phi(\beta)}-\varphi'(0+)W^{(\beta)}(a).
 \end{equation}
and
 \begin{equation}
  E_{-z}(\text{e}^{-\beta T^+_a}, T_a^+< T_{-b}^-)=\frac{W^{(\beta)}(b-z)}{W^{(\beta)}(a+b)}.
  \end{equation}
Moreover, using a fact in Kyprianou and Palmowski [26]   we get
 \begin{equation}
 P_{-z}^{(\Phi(\beta))}(T_{-b}<\infty)=1-\varphi'_{\Phi(\beta)}(0+)W_{\Phi(\beta)}(b-z).
 \end{equation}
Since (cf. Pistorius [32])
$$W^{(q)}(x)=e^{\Phi(q)x}W_{\Phi(q)}(x),\;\;\;
\varphi_{\Phi(q)}(\lambda)=\varphi(\Phi(q)+\lambda)-q,$$ we can
rewrite (3.15) as
 \begin{equation}
 P_{-z}^{(\Phi(\beta))}(T_{-b}<\infty)
=1-\varphi'(\Phi(\beta))\text{e}^{-\Phi(\beta)(b-z)}W^{(\beta)}(b-z).
\end{equation}
Inserting (3.12)-(3.16)  in (3.11) completes the proof.

The following result generalized the corresponding result  in Dos Reis [14] and Zhang and Wu [38] in which  the
classical compound Poisson risk model  and
the classical compound Poisson risk model perturbed by Brownian
motion were considered, respectively.  A different approach can be found in Landriault, Renaud and  Zhou [21].

\begin {theorem}  Suppose that the L\'evy process $X$, with
the Laplace exponent (1.1), drifts to $\infty$.    Let
$D=\int_0^{\infty}1(X(t)<0)dt$ denote the total duration for $X$ to
stay below 0.
Then  for $x>0$, $\beta>0$,
\begin{equation}
E_x\text{e}^{-\beta
D}=\varphi'(0+)\Phi(\beta)e^{\Phi(\beta)x}\int_x^{\infty}e^{-\Phi(\beta)y}W(y)dy.
\end{equation}
In particular,
$$E_0\text{e}^{-\beta
D}=\varphi'(0+)\frac{\Phi(\beta)}{\beta}.$$
\end {theorem}
{\bf Proof}\ The ideas of this proof were partly motivated by Dos Reis
[14] and Zhang and Wu [36].  For $\varepsilon\ge 0$, define,
with the convention that $\inf\emptyset=\infty$,
$$L_1(\varepsilon)=\inf\{t\ge 0:X(t)<-\varepsilon\},\;\;
R_1(\varepsilon)=\inf\{t\ge L_1(\varepsilon):X(t)=0\}.$$ In general,
for $k\ge 2$ recursively define
$$L_k(\varepsilon)=\inf\{t\ge R_{k-1}(\varepsilon):X(t)<-\varepsilon\},\;\;
R_k(\varepsilon)=\inf\{t\ge L_k(\varepsilon):X(t)=0\}.$$ If there
exists  some $k$ such that $\{t\ge
R_{k-1}(\varepsilon):X(t)<-\varepsilon\}=\emptyset$, then we define
$L_k(\varepsilon)=\infty$ (and consequently
$R_k(\varepsilon)=\infty)$ and $R_k-L_k=0$.

We first consider  the case where the paths of $X$ are of  bounded variation. For convenience we shall write  $L_k$ in place of $L_k(0)$ and  $R_k$ in place of $R_k(0)$. As the paths of
$X$ are of bounded variation, then $0\le T=L_1<R_1< L_2<R_2<\cdots,$
and $R_k-L_k$ represents the duration of the period of the surplus
from $k$-th below the level $0$ to the time that $X(t)$ first visits
at $0$ after $L_k$. Thus the random variable $D$ can be decomposed
as follows:
$$D=\sum_{k=1}^{N}(R_k-L_k),$$
where $N=\sup\{k:L_k<\infty\}$ ($N=0$ if the set is empty).
Note that $N$ has a geometric distribution,
\begin{equation}P_x(N=n)=
  \left\{
  \begin{array}{ll}
     R(x), &n=0,\\
    \psi(x)R(0)(\psi(0))^{n-1}, &n=1,2,\cdots,
  \end{array}
  \right.\nonumber
\end{equation}
where
$R(x)=1-\psi(x)$ and $\psi(x)=P(I_{\infty}<0|X(0)=x).$ The
stationarity and independence of increments of $X$ imply that given
$N=n$, $\{R_k-L_k, k=1,\cdots,n\}$ are mutually independent and
$\{R_k-L_k, k=2,\cdots,n\}$ are identically distributed.
A simple argument by using the law of double expectation  and the
strong Markov property yields
  \begin{eqnarray}  E_x(\text{e}^{-\beta D})
&=&P_x(N=0)E_x(\text{e}^{-\beta D}|N=0)\nonumber\\
&&+\sum_{n=1}^{\infty}P_x(N=n) E_x(\text{e}^{-\beta (R_1-L_1)}|N=n)\nonumber\\
&&\times E_x \left(e^{-\beta\sum_{k=2}^n
(R_k-L_k)}|N=n\right)\nonumber\\
&=&P_x(N=0)+\sum_{n=1}^{\infty}P_x(N=n) E_x(\text{e}^{-\beta
(R_1-L_1)}|L_n<\infty, L_{n+1}=\infty)\nonumber\\
&&\times \left\{E_x \left(e^{-\beta(R_2-L_2)}|L_n<\infty,L_{n+1}=\infty\right)\right\}^{n-1}\nonumber\\
&=&R(x)+R(0)E_x \left(e^{X(T)\Phi(\beta)}, T<\infty\right)
\sum_{n=1}^{\infty}\left\{E_0
\left(e^{X(T)\Phi(\beta)},T<\infty\right)\right\}^{n-1}\nonumber\\
&=&R(x)+E_x \left(e^{X(T)\Phi(\beta)}, T<\infty\right) \frac{R(0)
}{1-E_0\left( e^{X(T)\Phi(\beta)},T<\infty\right)},
\end{eqnarray}
where, in the third step, we have used
 \begin{eqnarray} E_x(\text{e}^{-\beta
(R_1-L_1)}|L_n<\infty, L_{n+1}=\infty)&=&E_x(E_{X(T)} e^{-\beta R_1}, T<\infty)/P_x(T<\infty)\nonumber\\
&=&E_x \left(e^{X(T)\Phi(\beta)}, T<\infty\right)/P_x(T<\infty), \nonumber
\end{eqnarray}
and
 \begin{eqnarray} E_x(\text{e}^{-\beta
(R_2-L_2)}|L_n<\infty, L_{n+1}=\infty)&=&E_0(E_{X(T)} e^{-\beta R_1}, T<\infty)/P_0(T<\infty)\nonumber\\
&=&E_0 \left(e^{X(T)\Phi(\beta)}, T<\infty\right)/P_0(T<\infty). \nonumber
\end{eqnarray}
By using
(2.9) one finds that
$$E_x\left( e^{X(T)\Phi(\beta)},T<\infty\right)=\beta e^{\Phi(\beta) x}\int_x^{\infty}e^{-\Phi(\beta)
y}W(y)dy-\frac{\beta}{\Phi(\beta)}W(x).$$ The result (3.17) follows,
since $R(x)=\varphi'(0+) W(x)$.

Next we  consider  the case where the paths of $X$ are of  unbounded variation. Note that, for
$\varepsilon>0$,
 $0\le L_1(\varepsilon)\le R_1(\varepsilon)\le
L_2(\varepsilon)\le R_2(\varepsilon)\le \cdots,$ and
$R_k(\varepsilon)-L_k(\varepsilon)$ represents the duration of the
period of the surplus from $k$-th below the level $-\varepsilon$ to
the time that $X(t)$ first visits at $0$ after $L_k(\varepsilon)$.
Let
$$D(\varepsilon)=\sum_{k=1}^{N(\varepsilon)}(R_k(\varepsilon)-L_k(\varepsilon)),$$
where $N(\varepsilon)=\sup\{k:L_k(\varepsilon)<\infty\}$
($N(\varepsilon)=0$ if the set is empty), which  has a geometric
distribution,
\begin{equation}P_x(N(\varepsilon)=n)=
  \left\{
  \begin{array}{ll}
     R(x+\varepsilon), &n=0,\\
    \psi(x+\varepsilon)R(\varepsilon)(\psi(\varepsilon))^{n-1}, &n=1,2,\cdots,
  \end{array}
  \right.\nonumber
\end{equation}
 where, as before,  $R(x)=1-\psi(x)$ and
$\psi(x)=P(I_{\infty}<0|X(0)=x).$  As above, given
$N(\varepsilon)=n$, $\{R_k(\varepsilon)-L_k(\varepsilon),
k=1,\cdots,n\}$ are mutually independent and
$\{R_k(\varepsilon)-L_k(\varepsilon), k=2,\cdots,n\}$ are
identically distributed.

Using the same argument  as above we have
\begin{equation}
E_x(\text{e}^{-\beta D(\varepsilon)})
=R(x+\varepsilon)+
\frac{R(\varepsilon)E_x\left(
e^{X(T_{-\varepsilon}^{-})\Phi(\beta)},T^{-}_{-\varepsilon}<\infty\right)
}{1-E_0 \left(e^{X(T_{-\varepsilon}^{-})\Phi(\beta)}, T_{-\varepsilon}^{-}<\infty\right)}.
 \end{equation}
It follows from (2.9) that
$$\lim_{\varepsilon\to
0}\frac{R(\varepsilon)}{1-E_0\left(e^{X(T_{-\varepsilon}^{-})\Phi(\beta)}, T_{-\varepsilon}^{-}<\infty\right)}
=\varphi'(0+)\frac{\Phi(\beta)}{\beta}.$$ From the
right continuity of the sample paths of $X(t)$, we have
$\lim_{\varepsilon\to 0}D(\varepsilon)=D$. Thus the result (3.17)
follows by letting $\varepsilon\to 0$ in (3.19) and using (2.9) and
$R(x)=\varphi'(0+) W(x)$. This ends the proof.

{\bf Remark 3.2} For a spectrally one-sided L\'evy process,
the double-integral transforms  of the duration of stay
inside/outside the interval $(0,B)$ $(B>0$ is a constant) before a
fixed time  have been obtained by Kadankov and Kadankova [22].
However, our result can not deduced by the known result above.

\section  {Applications to insurance risk theory}
\setcounter{equation}{0}

Spectrally negative L\'evy processes have been considered recently
in  Refs. [4, 8, 16, 28, 30, 33], among others,
in the context of insurance risk models. Motivated by applications
in option pricing and risk management, and inspired recent
developments in fluctuation theory for L\'evy processes, Biffis and
Kyprianou [8] and Biffis and Morales [9] defined  an extended
version of the Gerber and Shiu expected discounted penalty function
introduced by Gerber and Shiu [18]. In addition to the surplus
before ruin and the deficit at ruin, it includes the information on
 the last minimum of the surplus before ruin $I_T$, where
$T=\inf\{t>0:X(t)<0\}$ denoting the ruin time of $X$. The analysis
of the result is mainly based on the quintuple law in Doney and
Kyprianou [13].

Motivated by them, we now consider the other generalized version of
the Gerber-Shiu expected discounted penalty function:
\begin{equation}
\phi(x;q,w)=E_x(e^{-q T}w(X(T-),|X(T)|,
S_{T-}, I_{T-})1(T<\infty)),
\end{equation}
  where $x\ge 0$ is the initial
surplus, $q\ge 0$ can be interpreted as a force of interest,
$w:\Bbb{R}^4\rightarrow [0,\infty)$ is bounded measurable function.
Using Theorem 3.1 we get the following corollary:

\begin{corollary} Suppose that $X$  drifts to $\infty$,
$W^{(q)}|_{(0,\infty)}\in C^{2}(0,\infty)$. Then the function
defined in (4.1) can be written as
$$\phi(x;q,w)=\int_{[0,\infty)^4}w(y,z,a,b)\left\{K_x^{(q)}(y,z,a,b)
+1_{(\sigma\neq0)}\delta(y,z,a-x,b)K_6\right\} dydzdadb,$$
where $\delta$ is the multidimensional Dirac Delta function, and
$$ K_x^{(q)}(y,z,a,b)= 1(y\ge b, a\ge x, z>0, b>0)
\pi(-z-y)\sum_{i=1}^5 K_i,$$
$$K_1=\frac{{W^{(q)}}'(x-b){W^{(q)}}'(a-y)}{
W^{(q)}(a-b)},$$
$$K_2=-\frac{{W^{(q)}}'(x-b)W^{(q)}(a-y)
{W^{(q)}}'(a-b)}{{W^{(q)}}^2(a-b)},$$
$$K_3=-\frac{W^{(q)}(x-b){W^{(q)}}''(a-b)W^{(q)}(a-y)}{{W^{(q)}}^2(a-b)},$$
$$K_4=-\frac{W^{(q)}(x-b){W^{(q)}}'(a-b){W^{(q)}}'(a-y)}{{W^{(q)}}^2(a-b)},$$
$$K_5=\frac{2{W^{(q)}}'(a-b)W^{(q)}(a-y){W^{(q)}}'(a-b)}{{W^{(q)}}^3(a-b)},$$
$$K_6=Z^{(q)}(x)-\frac{q
W^{(q)}(x)}{\Phi(q)}
-\int_0^{\infty}u^{(q)}(x,y)(\Pi(-y)-\Pi(-\infty))dy.$$
\end{corollary}

To end this section we rewrite  the generalized Dickson's formula
for  the Cram\'er-Lundberg risk process (see Gerber and Shiu [18])
and for jump-diffusion process (see Zhang and Wang [38]) in a more
appealing form  in terms of the probabilities of ruin or the scale
functions.

Consider the jump-diffusion risk process:
 \begin{eqnarray}
 X(t)=x+ct-\sum_{j=1}^{N(t)}X_j +\sigma B(t), t \ge 0,
 \end{eqnarray}
where $x$ is the insurer's initial capital, $c$ is the premium rate,
$\{N(t), t\ge 0\}$ is a Poisson process with parameter $\lambda$ and
$\{X_k\}_{k\ge 1}$ are independent random variables with common
distribution $P=1-\overline{P}$, which has density $p$, mean $\mu$
and $P(0)=0$, $\{B(t), t\ge 0\}$ is a Brownian motion. Moreover, we
assume that $c>\lambda\mu$ and, $\{N(t), t\ge 0\}$, $\{X_k\}_{k\ge
1}$ and $\{B(t), t\ge 0\}$ are assumed to be independent. When
$\sigma=0$, (4.2) is called the Cram\'er-Lundberg risk process.
Those two processes correspond to the cases of spectrally negative
L\'evy processes with $\Pi\{(-\infty,0)\}<\infty$ and, with or
without Gaussian component. For details of risk theory, we refer the
readers to Asmussen [2] and Rolski et al. [34].

Obviously, $X$ is a spectrally negative  L\'evy process with
$E\text{e}^{\alpha (X(t)-x)}=\text{e}^{t\varphi(\alpha)},$ where
$\varphi$ is defined as
$$\varphi(\alpha)=c\alpha+\frac12\sigma^2\alpha^2+\lambda
\{\hat{p}(\alpha)-1\}.$$

The following  generalized Dickson's formula  for  the
Cram\'er-Lundberg risk process is due to Gerber and Shiu [18,
(6.5) and (6.6)]:
\begin{equation} f_{q}(y|x)=
  \left\{
  \begin{array}{ll}
  f_{q}(x|0)\frac{e^{\Phi(q) x}-\Psi (x)}{1-\Psi(0)}, &y>x\ge 0,\\
  f_{q}(y|0)\frac{e^{\Phi(q) x}\Psi(x-y)-\Psi(x)}{1-\Psi(0)}, & 0<y\le x,
  \end{array}
  \right.
\end{equation}
where $\Psi(x)=E_x(\text{e}^{-qT+\Phi(q)
X(T)}{1}(T<\infty))$, and
$$f_{q}(y|0)=\lambda c^{-1}e^{-\Phi(q) y}(1-P(y)),\;\;\;
\Psi(0)=\lambda c^{-1} \int_0^{\infty}x\text{e}^{-\Phi(q)
x}p(x)dx.$$ Using (1.2), we can write $\Psi(x)$ as
$\Psi(x)=\text{e}^{\Phi(q) x}P^{(\Phi(q))}_x(T<\infty).$
Furthermore, $c(1-\Psi(0))=E^{(\Phi(q))}X(1)=\varphi'(\Phi(q)).$
 As a result,  we can rewrite  the
generalized Dickson's formula (4.3)   as
\begin{equation} f_{q}(y|x)=
  \left\{
  \begin{array}{ll}
   \lambda \overline{P}(y)\frac{e^{\Phi(q)(x-y)}}{\varphi'(\Phi(q))}P^{(\Phi(q))}_x(T=\infty), & y>x\ge 0,\\
   \lambda\overline{P}(y)\frac{e^{\Phi(q)(x-y)}}{\varphi'(\Phi(q))}\left(P^{(\Phi(q))}_x(T=\infty)
   -P^{(\Phi(q))}_{x-y}(T=\infty)\right),
& 0<y\le x.
  \end{array}
  \right.
\end{equation}
The generalized Dickson's formula for jump-diffusion is due to Zhang
and Wang [39]:
\begin{equation} f_{q}(y|x)=
  \left\{
  \begin{array}{ll}
   \lambda \overline{P}(y)\Phi'(q)e^{-\Phi(q) y}\left(e^{\Phi(q)x}-M(x)\right),&  y>x\ge 0,\\
    \lambda \overline{P}(y)\Phi'(q)e^{-\Phi(q)
y}\left(e^{\Phi(q)x}M(x-y)-M(x)\right), & 0<y\le x,
  \end{array}
  \right.
\end{equation}
 where
$M(x)=E_x e^{-q T_0}.$ From Chiu and Yin [11, Theorem 2.3]   we
have
$$M(x)=E_x\left(e^{-q T+\Phi(q)
X(T)}1(T<\infty)\right).$$

By (1.2), $M(x)=\text{e}^{\Phi(q) x}P^{(\Phi(q))}_x(T<\infty).$
Therefore, (4.5) can be rewritten as
\begin{equation} f_{q}(y|x)=
  \left\{
  \begin{array}{ll}
\lambda \overline{P}(y)\frac{e^{\Phi(q)(x-y)}}{\varphi'(\Phi(q))}P^{(\Phi(q))}_x(T=\infty), & y>x\ge 0,\\
\lambda\overline{P}(y)\frac{e^{\Phi(q)(x-y)}}{\varphi'(\Phi(q))}\left(P^{(\Phi(q))}_x(T=\infty)
-P^{(\Phi(q))}_{x-y}(T=\infty)\right),
& 0<y\le x,
  \end{array}
  \right.
\end{equation}
where we have used
$\Phi'(q)\cdot\varphi'(\Phi(q))=1$.

{\bf Remark 4.1}  The results (4.4) and (4.6) can also be expressed
in terms of the scale functions by using (1.2). We  remark that the
positive safety loading condition is not required   in the case
$q>0$. Thus the corresponding conditions can be removed here.
\vskip0.3cm

\noindent\bf{\footnotesize Acknowledgements}\quad\rm
{\footnotesize Both authors thank two anonymous referees for their constructive suggestions
 which have led to much improvement on the paper.
 The first author is grateful to Professor Xiaowen Zhou  for useful discussions. The research of  Yuen was
supported by a university research grant of the University of Hong
Kong. The research of  Yin was supported by the National
Natural Science Foundation of China (No. 11171179),   the Research Fund for
the Doctoral Program of Higher Education of China (No. 20133705110002) and the Program for  Scientific Research Innovation Team in Colleges and Universities of Shandong Province.}\\[4mm]

\noindent{\bbb{References}}

\end{document}